\documentclass[preprint,reqno,12pt]{amsart}
\usepackage{amsfonts}

\usepackage{graphicx}
\usepackage[spanish,activeacute]{babel}
\usepackage[dvips]{}

\newtheorem{thm}{Theorem}[section]

\theoremstyle{definition}
\newtheorem{defn}[thm]{Definition}
\theoremstyle{remark}
\newtheorem{rem}[thm]{Remark}
\numberwithin{equation}{section}

\begin{document}
\title[]{New Trigonometric form  of The Hamilton's Quaternions}
\author{Mijail Andr'es Saralain Figueredo}
\address{Facultad de Matem$\acute{a}$tica, F$\acute{i}$sica y Computaci$\acute{o}$n,
Universidad Central Marta Abreu de Las Villas,Apartado Postal 54830,Santa
Clara ,Villa Clara,Cuba}
\email{mijail@mfc.uclv.edu.cu}
\thanks{mijail@mfc.uclv.edu.cu}
\date{\today}

\begin{abstract}
Is it possible to define, for certain values $n$ the product of vectors of
the real vector space of n dimensions , such that this is, with respect to
multiplication and the ordinary addition of vectors, a numerical system
which contains the system of real numbers? It can be proven that this cannot
be done. In the space of four dimensions this construction is possible if we
are apart from the commutativity of the multiplication. The resulting system
is the one of \textbf{QUATERNIONS}. In this work I first do a reminder of
the fundamental concepts of Hamilton's Hypercomplex and then a deep work
with such concepts.
\end{abstract}

\maketitle

%

%


\section{Define quaternion}

\begin{defn}
We shall call quaternions or simply Hamilton's hypercomplexes to an
expression in the form:
\begin{equation*}
Q=a+bi+cj+dk
\end{equation*}
where: $a,b,c,d\in \mathbb{R}$.Besides $i,j,k$ are imaginary units, pairwise
solutions of the equation $x^{2}=-1$, satisfying:
\begin{equation*}
ij=k=-ji
\end{equation*}
\begin{equation*}
jk=i=-kj
\end{equation*}
\begin{equation*}
ki=j=-ik
\end{equation*}
\begin{equation*}
ijk=-1
\end{equation*}
\begin{equation*}
i^{2}=j^{2}=k^{2}=-1.
\end{equation*}
\end{defn}

\begin{defn}
We shall say a quaternion is purely imaginary, if the first element of the
expression is equal to zero $(a=0,Im(\mathbb{H}))$.
\end{defn}

\begin{thm}
We say $Q=Q^{\prime }$, with $Q^{\prime }=a^{\prime }+b^{\prime }i+c^{\prime
}j+d^{\prime }k$ $(Q,Q^{\prime }\in \mathbb{H})$, ie two quaternions are
equal, if and only if, are equivalent the components of their imaginary and
real parts: $a=a^{\prime },b=b^{\prime },c=c^{\prime },d=d^{\prime }.$ \emph{%
\ }
\end{thm}

Proof. It's a direct consequence of the equality in $\mathbb{R}^{n}.$

\subsection{Fundamental Definition}

\begin{itemize}
\item  The sum and substraction are defined component by component, ie:
\begin{equation}
Q+Q^{\prime }=(a+a^{\prime })+(b+b^{\prime })i+(c+c^{\prime })j+(d+d^{\prime
})k
\end{equation}
\begin{equation}
Q-Q^{\prime }=(a-a^{\prime })+(b-b^{\prime })i+(c-c^{\prime })j+(d-d^{\prime
})k.
\end{equation}
\end{itemize}

\begin{itemize}
\item  The product is defined in the way:
\begin{equation}
Q^{\prime \prime }=(a+bi+cj+dk)(a^{\prime }+b^{\prime }i+c^{\prime
}j+d^{\prime }k)
\end{equation}
resulting:

$Q^{\prime\prime}=a a^{\prime}-b b^{\prime}- c c^{\prime}- d d^{\prime}+(a
b^{\prime}+a^{\prime}b + c d^{\prime}- c^{\prime}d) i + (a c^{\prime}+
a^{\prime}c -b d^{\prime}+ b^{\prime}d) j + (a d^{\prime}+ a^{\prime}d + b
c^{\prime}- b^{\prime}c) k$, with:
\begin{equation*}
a^{\prime\prime}= a a^{\prime}- b b^{\prime}- c c^{\prime}- d d^{\prime}
\end{equation*}
\begin{equation*}
b^{\prime\prime}= a b^{\prime}+ a^{\prime}b + c d^{\prime}- c^{\prime}d
\end{equation*}
\begin{equation*}
c^{\prime\prime}= a c^{\prime}+ a^{\prime}c - b d^{\prime}+ b^{\prime}d
\end{equation*}
\begin{equation*}
d^{\prime\prime}= a d^{\prime}+ a^{\prime}d + b c^{\prime}-b^{\prime}c
\end{equation*}
\begin{eqnarray}
Q^{\prime\prime}= a^{\prime\prime}+ b^{\prime\prime}i + c^{\prime\prime}j +
d^{\prime\prime}k.
\end{eqnarray}
\end{itemize}

\subsection{Fundamental Properties}

\begin{itemize}
\item  Commutativity of the sum: $Q_{1}+Q_{2}=Q_{2}+Q_{1}$.

\item  Associativity of the sum: $(Q_{1}+Q_{2})+Q_{3}=Q_{1}+(Q_{2}+Q_{3})$.

\item  Associativity of the product: $(Q_{1}\cdot Q_{2})\cdot
Q_{3}=Q_{1}\cdot (Q_{2}\cdot Q_{3})$.

\item  Distributivity: $(Q_{1}+Q_{2})\cdot Q_{3}=(Q_{1}+Q_{3})(Q_{2}+Q_{3})$.
\end{itemize}

\begin{rem}
Among the upper properties the following is missing:
\begin{equation}
Q_{1}Q_{2}\neq Q_{2}Q_{1}(in general).
\end{equation}
A very important result is Hamilton's hypercomplex are not a commutative
field.
\end{rem}

\begin{defn}
We shall call conjugate of a quaternion $Q$, and denote $\overline{Q}$, the
number:
\begin{equation}
\overline{Q}=a-bi-cj-dk.
\end{equation}
\end{defn}

Let us write now the sum and the difference of a quaternion with its
conjugate:

\begin{itemize}
\item  Sum:
\begin{equation*}
Q+\overline{Q}=(a+bi+cj+dk)+(a-bi-cj-dk)=2a=2\otimes \mathbb{R}e(\mathbb{H}).
\end{equation*}

\item  Difference:
\begin{equation*}
Q-\overline{Q}=(a+bi+cj+dk)-(a-bi-cj-dk)=2(bi+cj+dk)=2\otimes \mathbb{I}m(%
\mathbb{H}).
\end{equation*}
\end{itemize}

\begin{itemize}
\item  Product by its conjugate:
\begin{equation}
Q\cdot \overline{Q}=(a+bi+cj+dk)(a-bi-cj-dk)=a^{2}+b^{2}+c^{2}+d^{2}.
\end{equation}
\end{itemize}

\subsubsection{Other properties}

\begin{itemize}
\item  Selfpowered:
\begin{equation*}
\overline{\overline{Q}}=Q
\end{equation*}

\item  Additivity:
\begin{equation*}
\overline{Q_{1}+Q_{2}}=\overline{Q_{1}}+\overline{Q_{2}}
\end{equation*}

\item  Multiplicativity:
\begin{equation*}
\overline{Q_{1}\cdot Q_{2}}=\overline{Q_{2}}\cdot \overline{Q_{1}}
\end{equation*}

\item  Divisibility:
\begin{equation*}
\overline{(\frac{Q_{1}}{Q_{2}})}=\frac{\overline{Q_{1}}}{\overline{Q_{2}}}.
\end{equation*}
\end{itemize}

Is it possible to establish the inverse for the sum and the multiplication?

\begin{defn}
From $Q_{1}+Q_{2}=0$ we have: $a+a^{\prime }=0$, $b+b^{\prime }=0$, $%
c+c^{\prime }=0$, $d+d^{\prime }=0$. This implies $a^{\prime }=-a$, $%
b^{\prime }=-b$, $c^{\prime }=-c$, $d^{\prime }=-d$ . Therefor $%
Q_{2}=-a-bi-cj-dk$ . $Q_{2}=-Q_{1}.$
\end{defn}

Analogously,

\begin{defn}
If $Q_{1}\cdot Q_{2}=1$ with $Q_{1}\neq 0$ then $Q_{2}$ will be the
multiplicative inverse for $Q_{1}$.
\begin{equation}
Q_{2}=\frac{a-bi-cj-dk}{a^{2}+b^{2}+c^{2}+d^{2}}
\end{equation}
with $Q_{1}=a+bi+cj+dk$.
\end{defn}

Now, we've already defined the conjugate and multiplicative inverse. Let's
define division.

\begin{defn}
Division is define like this:
\begin{equation}
\frac{Q_{1}}{Q_{2}}=Q_{1}\cdot Q_{2}^{-1}
\end{equation}
\end{defn}

\begin{rem}
Notice when dividing we multiply the numerator by the multiplicative of the
denominator(recall (1.5))
\end{rem}

\section{Absolute Value}

\begin{defn}
We shall call absolute value or modulus of the quaternion, the nonnegative
real number $\left| Q\right| $,
\begin{equation*}
|Q|=|a+bi+cj+dk|=\sqrt{a^{2}+b^{2}+c^{2}+d^{2}}.
\end{equation*}
\end{defn}

Evidently, if we want to find the modulus of any quaternion: $Q+Q^{\prime
},Q^{\prime \prime },Q-Q^{\prime }$ ,etc; this will be the square root of
the sum of the squares of the real elements of each imaginary unit. Let us
notice that:
\begin{equation*}
|Q|^{2}=Q\cdot \overline{Q}.
\end{equation*}
As seen in formula $(1.7)$.

Thus, $|Q|^{2}=a^{2}+b^{2}+c^{2}+d^{2}$.

Now with $Q=a+bi+cj+dk$ and $\overline{Q}=a-bi-cj-dk$ (finding the conjugate
of the earlier):

\begin{itemize}
\item  $|\overline{Q}|=|Q|$

\item  $|Q_{1}\cdot Q_{2}|=|Q_{1}||Q_{2}|$

\item  ${|Q_{1}Q_{2}|}^{2}={|Q_{1}|}^{2}\cdot {|Q_{2}|}^{2}.$
\end{itemize}

\subsection{Norm}

\begin{defn}
The norm will be defined like this:
\begin{equation}
\|Q\|^{2}=Q\cdot \overline{Q}=|Q|^{2}
\end{equation}
\end{defn}

Now with $Q=a+bi+cj+dk$ and $\overline{Q}=a-bi-cj-dk$ (finding the conjugate
of the earlier):

\begin{itemize}
\item  $\|\overline{Q}\|=\|Q\|$

\item  $\|Q_{1}\cdot Q_{2}\|=\|Q_{1}\|\|Q_{2}\|^{2}$

\item  ${\|Q_{1}Q_{2}\|}^{2}={\|Q_{1}\|}^{2}\cdot {\|Q_{2}\|}.$
\end{itemize}

It is simple consequence of the modulus.

\begin{rem}
We defined division previously. Now we state the following equality:
\begin{equation*}
\frac{Q_{1}}{Q_{2}}=\frac{Q_{1}\cdot \overline{Q_{2}}}{\|Q_{2}\|}
\end{equation*}
\end{rem}

\begin{defn}
We shall call unit quaternion the Hamilton's hypercomplex which satisfies:
\begin{equation}
\|Q\|=1
\end{equation}
\end{defn}

i.e. $a^{2}+b^{2}+c^{2}+d^{2}=1$.

\section{Several ways of defining a quaternion}

\subsection{Vector Form}

\begin{center}
$Y: \mathbb{H} \longrightarrow {\mathbb{R}}^{4}$

$a + b i + c j + d k \longrightarrow (a, b, c, d);a, b, c, d \in\mathbb{R}$
\end{center}

Let us prove $Y$ is a bijective application.

$\forall x_{1},x_{2}\in \mathbb{H}$ with $x_{1}\neq x_{2}\Rightarrow
Y(x_{1})\neq Y(x_{2})$.

- Therefor is injective.

$\forall y \in \mathbb{R^{4}} \exists x \in \mathbb{H}:Y(x)=y$.

- Therefor is surjective.

We've proved $Y$ is a bijective function.

Provided that$(\mathbb{H},+)$ is an Abelian group, with the sum as defined
earlier; that the product is distributive with respect to the sum, $%
Q_{1}(Q_{2}+Q_{3})=Q_{1}\cdot Q_{2}+Q_{1}\cdot Q_{3}$ and the product is
associative $Q_{1}(Q_{2}Q_{3})=(Q_{1}Q_{2})Q_{3}$; then, quaternions with
the operations of sum and product define a ring.

$(\mathbb{H},+,\ast )$ is a ring.

\subsubsection{Definitions of this notation}

\begin{enumerate}
\item  The quaternion$(0,0,0,0)$ is the neutral quaternion for the sum. It's
obvious from the definition of sum: $(a,b,c,d)+(0,0,0,0)=(a,b,c,d)$.

\item  Analogously, the neuter for the product is $(1,0,0,0)$ because $%
(a,b,c,d)(1,0,0,0)=(a,b,c,d).$

\item  Noting we are using our defined product: $\alpha (a,b,c,d)=(\alpha
a,\alpha b,\alpha c,\alpha d).$

\item  $(0,1,0,0)(0,1,0,0)=(-1,0,0,0)$ as we can see this quaternion is
identified by $i^{2}$.
\end{enumerate}

From now on, we'll consider: $Y(1)=(1,0,0,0)$ $Y(i)=(0,1,0,0)$ $%
Y(j)=(0,0,1,0)$ $Y(k)=(0,0,0,1)$.

We wonder,is it an isomorphism? We only have to prove the following holds:

$Y [( a + b i + c j + d k)+( s + r i + t j + h k)] = Y [a + b i + c j + d
k]+Y [ s + r i + tj + h k]$.

\subsection{Other Vector Form}

\begin{center}
$\Gamma: \mathbb{H} \longrightarrow {\mathbb{R}}^{4}$

$a + b i + c j + d k \longrightarrow (a,\overrightarrow{v});a\in\mathbb{R},%
\overrightarrow{v}$ vector in $\mathbb{R}^{3}$
\end{center}

The conjugate of this vector is: $\overline{(a,\overrightarrow{v})}=(a,-%
\overrightarrow{v})$. I remark the multiplication with this notation is: $%
Q_{1}\cdot Q_{2}=(a,\overrightarrow{v})(a^{\prime},\overrightarrow{v_{1}}%
)=(aa^{\prime}-\overrightarrow{v}\cdot \overrightarrow{v_{1}},a%
\overrightarrow{v_{1}}+a^{\prime}\cdot \overrightarrow{v}+ \overrightarrow{v}
\times \overrightarrow{v_{1}})$

\subsection{Matrix form}

The matrix form of defining a quaternion is:

\begin{center}
$\Omega : \mathbb{H} \longrightarrow \mathbb{M}$

$a + b i + c j + d k \longrightarrow \left(
\begin{array}{cc}
a+bi & c+di \\
-c+di & a-bi
\end{array}
\right)$
\end{center}

In order to illustrate this notation, it's convenient to develop the
following, by taking
\begin{equation*}
\Omega (1)=\left(
\begin{array}{cc}
1 & 0 \\
0 & 1
\end{array}
\right) ,\Omega (i)=\left(
\begin{array}{cc}
i & 0 \\
0 & -i
\end{array}
\right) ,\Omega (j)=\left(
\begin{array}{cc}
0 & 1 \\
-1 & 0
\end{array}
\right) ,\Omega (k)=\left(
\begin{array}{cc}
0 & i \\
i & 0
\end{array}
\right) .
\end{equation*}
These are called Pauli's matrixes.

The quaternion should be then:
\begin{equation*}
\Omega (a+bi+cj+dk)=\left(
\begin{array}{cc}
a+bi & c+di \\
-c+di & a-bi
\end{array}
\right)
\end{equation*}
\begin{equation*}
\Omega (a+bi+cj+dk)=a\cdot \left(
\begin{array}{cc}
1 & 0 \\
0 & 1
\end{array}
\right) +b\cdot \left(
\begin{array}{cc}
i & 0 \\
0 & -i
\end{array}
\right) +c\cdot \left(
\begin{array}{cc}
0 & 1 \\
-1 & 0
\end{array}
\right) +d\cdot \left(
\begin{array}{cc}
0 & i \\
i & 0
\end{array}
\right)
\end{equation*}
\begin{equation*}
\Omega (a+bi+cj+dk)=a \Omega (1)+b \Omega (i)+c \Omega (j)+d \Omega (k).
\end{equation*}
We have only left to wonder, is this function an isomorphism? It is, indeed,
and the proof is equivalent to the earlier function $(Y)$.If we find the
determinant of the matrix $\left(
\begin{array}{cc}
a+bi & c+di \\
-c+di & a-bi
\end{array}
\right) $ and calculate the modulus of the quaternion $a+bi+cj+dk$ we can
see identical results, therefor it is easily seen this function is an
isomorphism.

\subsection{Trigonometric form}

\begin{center}
$\Pi :\mathbb{H}\longrightarrow \mathbb{T}$ ,$\mathbb{T}:$Trigonometric form.

$a + b i + c j + d k \longrightarrow \rho·Cis\theta + (\rho_{0}Cis\beta) j $.
\end{center}

To illustrate this definition we write:

$a+bi=\rho Cis\theta $,$c+di=\rho _{0}Cis\beta $,$-c+di=-\rho _{0}Cis(-\beta
)$,$a-bi=\rho Cis(-\theta )$ with $a+bi,c+di,-c+di,a-bi\in \mathbb{C}.$

Finding the modulus of the earlier complex number we obtain:

$\rho = | a + b i | = | a - b i | =\sqrt{a^{2}+b^{2}} , \rho_{0} = | c +di |
= | -c + di | =\sqrt{c^{2}+d^{2}}.$

\subsection{New Trigonometric form}

\begin{center}
$\Upsilon :\mathbb{H}\longrightarrow \mathbb{J}$ ,$\mathbb{J}:$New
trigonometric form.

$a + b i + c j + d k \longrightarrow \rho(Cos(\alpha)+Sin(\alpha)j)$.

$a + b i + c j + d k \longrightarrow \rho·Cjs(\alpha) $.
\end{center}

\begin{rem}
The following is the short way to express the earlier. We know the field of
the quaternions is not commutative, so all of the transformations are
equivalent.
\begin{equation}
\rho·Cjs(\alpha )\equiv \rho (Cos(\alpha )+Sin(\alpha ))j
\end{equation}
The idea of the development is to shorten the trigonometric form.
\end{rem}

From complex analysis we know that:
\begin{equation*}
Ln(z)=Ln|z|+iArg(z)
\end{equation*}
\begin{equation*}
Tan^{-1}(z)=\frac{i}{2}Ln(\frac{i+z}{i-z})
\end{equation*}

\begin{defn}
Applying Pythagorean Theorem in $\mathbb{R}^{4}$, ie, working with the
complex axes(Fig.1), we define:
\begin{equation}
\rho =\sqrt{(a+bi)^{2}+(c+di)^{2}}
\end{equation}
\end{defn}

\begin{equation}
a+bi=\rho Cos(\alpha );c+di=\rho Sen(\alpha );Tan(\alpha )=\frac{c+di}{a+bi}
\end{equation}
From (3.3) we have that:
\begin{equation*}
\alpha =Tan^{-1}(\frac{c+di}{a+bi})=Cos^{-1}(\frac{a+bi}{\rho })=Sin^{-1}(%
\frac{c+di}{\rho })
\end{equation*}
Now, for obtaining the true value of $(\alpha )$ any of the three following
trigonometric function value's must be calculated: $%
Tan^{-1},Cos^{-1},Sin^{-1}$.They must be equal to $(\alpha )$.Now we wonder,
why not to consider the value of $(\rho )$ as the modulus of the quaternion?
The answer is the following: The trigonometric equality$(3.10)$ does not
always hold, only when: $b=d=0$. And the hypercomplex were not defined under
these conditions. The conclusion I reach about it, is that if these
equalities hold then we can take the modulus of the quaternion as $(\rho )$.
From (3.7) we have that:
\begin{equation*}
Tan^{-1}(\frac{c+di}{a+bi})=\frac{i}{2}Ln(\frac{a+d+(b-c)i}{a-d+(b+c)i})
\end{equation*}
ie:
\begin{equation*}
\alpha =\frac{i}{2}Ln(\frac{a+d+(b-c)i}{a-d+(b+c)i})
\end{equation*}
As we already can see we have $\rho $ which gives us (3.8), $(\alpha )$
which yields us (3.10) or (3.12). Then we can write a hypercomplex as:
\begin{equation*}
a+bi+cj+dk=\rho Cjs(\alpha )
\end{equation*}
where $\rho ,\alpha \in \mathbb{C}.$

\begin{center}
\includegraphics[width=10cm,height=7cm]{img.bmp}
\end{center}

In the earlier graph we can see that a quaternion can be represented in a
system of two complex coordinates. It is impossible to represent the set $%
\mathbb{C}$ in a coordinate axis. I only do it to see where I obtain $(3.8)$
and $(3.9)$. At the same time, each complex axis I show has values in $%
\mathbb{R}^{2}$,ie, we have the earlier graph in $\mathbb{R}^{4}.$

The remark I will write further is just a question. When defining a field in
trigonometric form, do we have to reference the before field?

As we can see I've arrived to a new way of defining quaternion. Now we have:
\begin{equation}
\rho Cjs(\alpha )\in \mathbb{H}
\end{equation}
with $\rho ,\alpha \in \mathbb{C}$, then we can write $\rho =\rho
_{1}(Cos(\beta )+iSin(\beta ))$ with $\rho _{1}\in \mathbb{R},\alpha \in
\mathbb{C}$. So
\begin{equation*}
a+bi+cj+dk=\rho Cjs(\alpha )=\rho _{1}(Cos(\beta )+iSin(\beta ))(Cos(\alpha
)+Sin(\alpha )j)
\end{equation*}
From complex analysis we know that:
\begin{equation*}
Sen(\alpha )=Sen(x+iy)=Sen(x)Cosh(y)+iCos(x)Senh(y)
\end{equation*}
\begin{equation*}
Cos(\alpha )=Cos(x+iy)=Cos(x)Cosh(y)-iSen(x)Senh(y)
\end{equation*}
Substituting

$a+bi+cj+dk=\rho _{1}[Cos(\beta )+iSin(\beta
)][Cos(x)Cosh(y)-iSen(x)Senh(y)+[Cos(x)Cosh(y)-iSen(x)Senh(y)]j].$

yields:
\begin{equation*}
a=\rho _{1}Cos(x)Cosh(y)
\end{equation*}
\begin{equation*}
b=-\rho _{1}Sin(x)Sinh(y)
\end{equation*}
\begin{equation*}
c=\rho _{1}Sin(x)Cosh(y)
\end{equation*}
\begin{equation*}
d=\rho _{1}Cos(x)Sinh(y)
\end{equation*}
William Rowan Hamilton in ''On a New Species of Imaginary Quantities
Connected with a Theory of Quaternions'', wrote:
\begin{equation*}
a=\rho Cos(\theta )
\end{equation*}
\begin{equation*}
b=\rho Sin(\theta )Cos(\vartheta )
\end{equation*}
\begin{equation*}
c=\rho Sin(\theta )Sin(\theta )Cos(\psi )
\end{equation*}
\begin{equation*}
d=\rho Sin(\theta )Sin(\theta )Sin(\psi )
\end{equation*}

\subsection{Trigonometric Matrix Form}

\begin{center}
$\Psi: \mathbb{H} \longrightarrow \mathbb{M}$

$a + bi+ cj + dk \longrightarrow \left(
\begin{array}{cc}
\rho Cis\theta & \rho_{0}Cis\beta \\
-\rho_{0} Cis\beta & \rho Cis(-\theta)
\end{array}
\right) $
\end{center}

I develop this kind of definition mostly for a practical work, such as
multiplying or dividing, due to these operations are easily performed in
trigonometric form.

Let's see that:

\begin{center}
$\left(
\begin{array}{cc}
\rho Cis\theta & \rho_{0} Cis\beta \\
-\rho_{0}Cis\beta & \rho Cis(-\theta)
\end{array}
\right)$ = $\left(
\begin{array}{cc}
\rho Cis\theta & 0 \\
0 & \rho Cis(-\theta)
\end{array}
\right)$ $+$ $\left(
\begin{array}{cc}
0 & \rho_{0} Cis\beta \\
-\rho_{0} Cis\beta & 0
\end{array}
\right) $.

$\left(
\begin{array}{cc}
\rho Cis\theta & \rho_{0} Cis\beta \\
-\rho_{0}Cis\beta & \rho Cis(-\theta)
\end{array}
\right)$ = $\left(
\begin{array}{cc}
\rho Cos\theta & 0 \\
0 & \rho Cos(-\theta)
\end{array}
\right) $ $+$ $\left(
\begin{array}{cc}
i\rho Sen\theta & 0 \\
0 & i\rho Sen(-\theta)
\end{array}
\right) $ $+$ $\left(
\begin{array}{cc}
0 & \rho_{0} Cos\beta \\
-\rho_{0} Cos\beta & 0
\end{array}
\right) $ $+$ $\left(
\begin{array}{cc}
0 & \rho_{0}Sen\beta \\
-\rho_{0}Sen\beta & 0
\end{array}
\right) $.
\end{center}

\subsection{Logarithmic Form}

Let us define the following map:

\begin{center}
$\Gamma: \mathbb{H} \longrightarrow \mathbb{L}$

$a+bi+cj+dk \longrightarrow \ln ab^{i}c^{j}d^{k}$
\end{center}

This way of definition is to convert the elements $a,b,c,d\in \mathbb{R}$ in
natural logarithms.
\begin{equation*}
a=\ln a,b=\ln b,c=\ln c,d=\ln d
\end{equation*}
(the upper is an abuse of notation) Obviously it's possible to think each
part of the quaternion as real logarithms, because they are real. This
notation is only possible when the real elements are not zero.

\subsection{Exponential Form}

\begin{center}
$\Delta: \mathbb{H} \longrightarrow \mathbb{E}$

$a+bi+cj+dk \longrightarrow \rho e^{i\theta}+\rho_{0}e^{i\beta}j.$
\end{center}

Let us define the exponential form . First, we know from complex that:$%
e^{z}=e^{x}(Cosy+iSeny)$,this means that
\begin{equation}
e^{yi}=Cosy+iSeny=Cisy.
\end{equation}

Then: $Cis\theta =e^{i\theta},Cis\beta =e^{i\beta}.$

As we know a $a+bi+cj+dk$ can be written as $\rho Cis\theta +\rho
_{0}Cis\beta j.$

$a+bi+cj+dk=\rho Cis\theta + \rho_{0} Cis\beta j$, by $(3.5),$ $%
a+bi+cj+dk=\rho Cis\theta + \rho_{0} Cis\beta j=\rho
e^{i\theta}+\rho_{0}e^{i\beta}j.$

Now calculate$(a+bi+cj+dk)^{n}.$

$(a+bi+cj+dk)^{n}=(\rho e^{i\theta }+\rho _{0}e^{i\beta }j)^{n},$developing
Newton's binomial we have: $(\rho e^{i\theta}+\rho_{0}e^{i\beta}j)^{n}=$

$=(\rho e^{i\theta})^{n}+ C_{n}^{1}(\rho
e^{i\theta})^{n-1}\rho_{0}e^{i\beta}j+C_{n}^{2}(\rho
e^{i\theta})^{n-2}(\rho_{0}e^{i\beta}j)^{2}+...C_{n}^{n-1}\rho
e^{i\theta}(\rho_{0}e^{i\beta}j)^{n-1}+(\rho_{0}e^{i\beta}j)^{n}. $
\begin{eqnarray}
(\rho e^{i\theta}+\rho_{0}e^{i\beta}j)^{n}= \sum_{h=0}^{n}\left(
\begin{array}{c}
n \\
h
\end{array}
\right)(\rho e^{i\theta})^{n-h}(\rho_{0}e^{i\beta}j)^{h}
\end{eqnarray}

\begin{eqnarray}
(\rho e^{i\theta}+\rho_{0}e^{i\beta}j)^{n}=(\rho
e^{i\theta})^{n}\sum_{h=0}^{n}\left(
\begin{array}{c}
n \\
h
\end{array}
\right)(\frac{\rho_{0}}{\rho}e^{i(\beta-\theta)}j)^{h}
\end{eqnarray}

Now we must verify by induction if $(3.7)$ holds.

Obviously it holds in the case $n=1.$ Suppose it holds for $n=k$ and prove
it is true in the case $n=k+1.$ Now, with the new definition in
trigonometric form we'll prove the following:

\begin{thm}
We know $a+bi+cj+dk=\rho Cjs(\alpha)$, then $(a+bi+cj+dk)^{n}=(\rho
Cjs(\alpha))^{n}=\rho^{n}Cjs(n \alpha).$
\end{thm}

\section{Functions of a hypercomplex variable}

If a variable $w$ is related with $z$ such that to each value of $z$ in $%
\mathbb{H}$ corresponds a value or set of defined values of $w$, then $w$ is
a function of the hypercomplex variable $z$, $w=f(z)$. If $z=a+bi+cj+dk$ and
$w=u+vi+sj+tk$ with the values of $a,b,c,d,u,v,s,t\in \mathbb{R}$. $%
u+vi+sj+tk=f(a+bi+cj+dk)$, and each of the real variables $u,v,s,t\in
\mathbb{R}$ are determined by the real quartet $a,b,c,d$. That is to say, $%
u=u(a,b,c,d)$,$v=v(a,b,c,d)$, $s=s(a,b,c,d)$,$t=t(a,b,c,d)$.

Example 1: $w=z^{2}+5.$

$u + v i + s j + t k = (a + b i + c j + d k)^{2} + 5,$

$u+vi+sj+tk=a^{2}-b^{2}-c^{2}-d^{2}+2abi+2acj+2adk+5.$ Then:

$u(a, b, c, d) = a^{2} - b^{2} - c^{2} - d^{2}+5,$

$v (a, b, c, d) = 2 a$, $s (a, b, c, d) = 2 a c$, $t(a, b, c, d) =2ad.$

\subsection{Limit Definition}

Let $f(z)$ a function defined in all the points in some neighborhood of $%
z_{0}$. We say that $w_{0}$ is the limit $f(z)$, when $z$ tends to $z_{0}$,
\begin{equation*}
\lim_{z\rightarrow z_{0}}f(z)=w_{0}.
\end{equation*}
That is, for all positive epsilon exists a positive number lambda such that:
\begin{equation*}
|f(z)-w_{0}|<\epsilon
\end{equation*}
when:
\begin{equation*}
|z-z_{0}|<\lambda (z\neq z_{0}).
\end{equation*}

Suppose that,
\begin{equation*}
\lim_{z\rightarrow z_{0}}f(z)=u_{0}+v_{0}i+s_{0}j+t_{0}k
\end{equation*}
where $f(z)=u+vi+sj+tk,$ $z=a+bi+cj+dk,$ $z_{0}=a_{0}+b_{0}i+c_{0}j+d_{0}k.$
Then by the inequality it becomes in:

\begin{equation*}
|u+vi+sj+tk-(u_{0}+v_{0}i+s_{0}j+t_{0}k)|<\epsilon
\end{equation*}
when:
\begin{equation*}
|a+bi+cj+dk-(a_{0}+b_{0}i+c_{0}j+d_{0}k)|<\lambda .
\end{equation*}
Making algebraic transformations, $(4.5),(4.6)$ yield:
\begin{equation*}
|(u-u_{0})+(v-v_{0})i+(s-s_{0})j+(t-t_{0})k|<\epsilon
\end{equation*}
when:
\begin{equation*}
|(a-a_{0})+(b-b_{0})i+(c-c_{0})j+(d-d_{0})k|<\lambda .
\end{equation*}
But from $(4.7)$ and $(4.8)$ we have, applying $(2.1)$:
\begin{equation*}
\sqrt{(u-u_{0})^{2}+(v-v_{0})^{2}+(s-s_{0})^{2}+(t-t_{0})^{2}}<\epsilon
\end{equation*}
when:
\begin{equation*}
\sqrt{(a-a_{0})^{2}+(b-b_{0})^{2}+(c-c_{0})^{2}+(d-d_{0})^{2}}<\lambda .
\end{equation*}

Uniqueness of the limit of a hypercomplex function: Suppose there exist two
limit points $w_{0},w_{1}$ $(w_{0}\neq w_{1}).$ By definition of limit:
\begin{equation*}
|f(z)-w_{0}|<\frac{\epsilon}{2} , when |z-z_{0}|<\lambda
\end{equation*}
\begin{equation*}
|f(z)-w_{1}|<\frac{\epsilon}{2} , when |z-z_{1}|<\lambda
\end{equation*}
Let's work with the real number$|w_{1}-w_{0}|.$ The following step is
obvious:

\begin{equation*}
|w_{1}-w_{0}|=|w_{1}-w_{0}+f(x)-f(x)|=|-f(x)+w_{1}+f(x)-w_{0}|=|-(f(x)-w_{1})+f(x)-w_{0}|\leq
\end{equation*}
\begin{equation*}
\leq|-(f(x)-w_{1})|+|f(x)-w_{0}|=|f(x)-w_{1}|+|f(x)-w_{0}|<\frac{\varepsilon%
}{2}+\frac{\varepsilon}{2}<\epsilon
\end{equation*}
\begin{equation*}
|w_{1}-w_{0}|<\epsilon ,
\end{equation*}
but as $w_{0},w_{1}$ are constant $|w_{1}-w_{0}|$ cannot be made as small as
one wants. Then, $w_{0}=w_{1}.$

Define a function $f(z)=w=u+vi+sj+tk$ such that:

\begin{equation*}
f^{\prime}( z ) = \lim_{\nabla a \rightarrow 0}\frac{u(a+\nabla a)-u(a)}{%
\nabla a}+\lim_{\nabla a \rightarrow 0}\frac{v(a+\nabla a)-v(a)}{\nabla a} i
+ \lim_{\nabla a \rightarrow 0}\frac{s(a+\nabla a)-s(a)}{\nabla a}j +
\lim_{\nabla a \rightarrow 0}\frac{t(a+\nabla a)-t(a)}{\nabla a}k.
\end{equation*}

\begin{equation*}
f^{\prime}( z ) = \lim_{\nabla b \rightarrow 0}\frac{u(b+\nabla b)-u(b)}{%
\nabla bi}+ \lim_{\nabla b \rightarrow 0}\frac{v(b+\nabla b)-v(b)}{\nabla bi}%
i + \lim_{\nabla b \rightarrow 0}\frac{s(b+\nabla b)-s(b)}{\nabla bi}j +
\lim_{\nabla b \rightarrow 0}\frac{t(b+\nabla b)-t(b)}{\nabla bi} k.
\end{equation*}

\begin{equation*}
f^{\prime}( z ) =\lim_{\nabla c \rightarrow 0}\frac{u(c+\nabla c)-u(c)}{%
\nabla cj} + \lim_{\nabla c \rightarrow 0}\frac{v(c+\nabla c)-v(c)}{\nabla cj%
} i+ \lim_{\nabla c \rightarrow 0}\frac{s(c+\nabla c)-s(c)}{\nabla cj}j +
\lim_{\nabla c \rightarrow 0}\frac{t(c+\nabla c)-t(c)}{\nabla cj}k.
\end{equation*}

\begin{equation*}
f^{\prime}( z ) =\lim_{\nabla d \rightarrow 0}\frac{u(d+\nabla d)-u(d)}{%
\nabla dk} + \lim_{\nabla d \rightarrow 0}\frac{v(d+\nabla d)-v(d)}{\nabla dk%
}i + \lim_{\nabla d \rightarrow 0}\frac{s(d+\nabla d)-s(d)}{\nabla dk}j +
\lim_{\nabla d \rightarrow 0}\frac{t(d+\nabla d)-t(d)}{\nabla dk} k.
\end{equation*}

Then:

\begin{equation}
f^{\prime}( z ) = \frac{\partial u}{\partial a}+ \frac{\partial v}{\partial a%
}i + \frac{\partial s}{\partial a}j + \frac{\partial t}{\partial a}k.
\end{equation}

\begin{equation}
f^{\prime}( z ) = - \frac{\partial u}{\partial b}i+ \frac{\partial v}{%
\partial b}+ \frac{\partial s}{\partial b}k -\frac{\partial t}{\partial b} j.
\end{equation}

\begin{equation}
f^{\prime}( z ) = - \frac{\partial u}{\partial c} j -\frac{\partial v}{%
\partial c} k +\frac{\partial s}{\partial c} + \frac{\partial t}{\partial c}
i.
\end{equation}

\begin{equation}
f^{\prime}( z ) = - \frac{\partial u}{\partial d} k +\frac{\partial v}{%
\partial d} j - \frac{\partial s}{\partial d}i + \frac{\partial t}{\partial d%
}.
\end{equation}

From the upper unqualified we conclude that a hypercomplex function is
analytic or integer if:

\begin{eqnarray}
\frac{\partial u}{\partial a}=\frac{\partial v}{\partial b}= \frac{\partial s%
}{\partial c}=\frac{\partial t}{\partial d}.
\end{eqnarray}

\begin{eqnarray}
\frac{\partial v}{\partial a} = -\frac{\partial u}{\partial b}=\frac{%
\partial t}{\partial c}=-\frac{\partial s}{\partial d}.
\end{eqnarray}

\begin{eqnarray}
\frac{\partial s}{\partial a}=-\frac{\partial t}{\partial b} = - \frac{%
\partial u}{\partial c} =\frac{\partial v}{\partial d}.
\end{eqnarray}

\begin{eqnarray}
\frac{\partial t}{\partial a}=\frac{\partial s}{\partial b} =-\frac{\partial
v}{\partial c} = - \frac{\partial u}{\partial d}.
\end{eqnarray}

We call the latter the generalized theorem of Cauchy-Riemann.

References:

[1] Kurochov, Algebra Lineal.

[2] Teresita Noriega, Algebra Lineal II.

[3] Ebbinghaus et al, Number, Stringer, 1991.

[4] Peter Grogono, Rotation with Quaternions, December 2001.


\end{document}